\documentclass[twoside]{article}
\usepackage{amssymb}
\usepackage{indentfirst}
\usepackage{dsfont}
\usepackage{ntheorem}
\usepackage{amsmath}

\begin{document}
\baselineskip 13pt \noindent \thispagestyle{empty}

\markboth{\centerline{\rm Z. Tu \; \& \; L. Wang}}{\centerline{\rm
Rigidity of proper holomorphic mappings}}

\begin{center} \Large{\bf  Rigidity of proper holomorphic mappings between \\
\vskip 4pt certain unbounded non-hyperbolic domains }
\end{center}

\begin{center}
\noindent\text{Zhenhan TU\; and \; Lei WANG$^{*}$}\\
\vskip 4pt \noindent\small {School of Mathematics and Statistics,
Wuhan
University, Wuhan, Hubei 430072, P.R. China} \\
\noindent\text{Email: zhhtu.math@whu.edu.cn (Z. Tu),\;}
{wanglei2012@whu.edu.cn (L. Wang)}
\renewcommand{\thefootnote}{{}}
\footnote{\hskip -16pt {$^{*}$Corresponding author. \\ } }
\end{center}

\begin{center}
\begin{minipage}{13cm}
{\bf Abstract.} {\small The Fock-Bargmann-Hartogs domain
$D_{n,m}(\mu)$ ($\mu>0$)
 in $\mathbf{C}^{n+m}$ is defined by the inequality $\|w\|^2<e^{-\mu\|z\|^2},$
 where $(z,w)\in \mathbf{C}^n\times \mathbf{C}^m$, which is an unbounded non-hyperbolic domain
 in $\mathbf{C}^{n+m}$. Recently, Yamamori gave an explicit formula for the Bergman
kernel of the Fock-Bargmann-Hartogs domains in terms of the
polylogarithm functions and Kim-Ninh-Yamamori determined the
automorphism group of the domain $D_{n,m}(\mu)$. In this article, we
obtain rigidity results on proper
 holomorphic mappings between two equidimensional Fock-Bargmann-Hartogs domains.
Our rigidity result implies that any proper holomorphic self-mapping
on the Fock-Bargmann-Hartogs domain $D_{n,m}(\mu)$  with $m\geq 2$
must be an automorphism.

 }
\end{minipage}
\end{center}

\section{Introduction}
 In 1977, Alexander \cite{Ale} proved the following fundamental result.
\newtheorem{lem}{Lemma}[section]
\newtheorem{thm}{Theorem}[section]
\newtheorem{prop}{Proposition}[section]

\vskip 6pt \noindent {\bf Theorem 1.A} (Alexander \cite{Ale}) {\it
If $f:\mathbf{B}^n\rightarrow \mathbf{B}^n$ ($n\geq 2$) is a proper
holomorphic self-mapping of the unit ball in $\mathbf{C}^n$, then
$f$ is an automorphism of $\mathbf{B}^n$.}

\vskip 6pt

Alexander's theorem has been generalized to several classes of
domains. Especially, there are many important results concerning
proper holomorphic mapping $f: D_1\rightarrow D_2$ between two
bounded pseudoconvex domains $D_1,\; D_2$ in $\mathbb{C}^n$ with
smooth boundary. If the proper holomorphic mapping $f$ extends
smoothly to the closure of $D_1$, then the extended mapping takes
the boundary $bD_1$ into the boundary $bD_2$, and it satisfies the
tangential Cauchy-Riemann equations on $bD_1$. Thus the proper
holomorphic mapping $f:D_1\rightarrow D_2$ leads naturally to the
geometric study of the mappings from $bD_1$ into $bD_2$. These
researches are often heavily based on analytic techniques about the
mapping on boundaries (e.g., see Forstneri\v{c} \cite{F} and Huang
 \cite{H} for references). In this regard, respectively, Diederich and
Forn{\ae}ss \cite{Diederich} and Bedford and Bell \cite{Bedf} proved
the following results.

\vskip 6pt \noindent {\bf Theorem 1.B} (Diederich and Fornaess
\cite{Diederich}) {\it If $\Omega, D\subset \mathbf{C}^n(n\geq 2)$
are smoothly bounded pseudoconvex domains and $\Omega$ is strongly
pseudoconvex, then any proper holomorphic mapping $f$ of $\Omega$
into $D$ is a local biholomorphism. Thus, if $D$ is simply
connected, then the mapping $f$ is biholomorphic. }

\vskip 6pt \noindent {\bf Theorem 1.C} (Bedford and Bell
\cite{Bedf}) {\it Let $D$ be bounded weakly pseudoconvex domain in
$\mathbf{C}^n(n\geq 2)$ with smooth real-analytic boundary. Then any
proper holomorphic self-mapping of $D$ is an automorphism. }

\vskip 6pt

We remark that $f(z_1,z_2)=(z_1,z_2^2):\; |z_1|^2+|z_2|^4<1
\rightarrow |w_1|^2+|w_2|^2<1$ is a proper holomorphic mapping
between two bounded pseudoconvex domains in $\mathbb{C}^2$ with
smooth real-analytic boundary, but it is  branched and is not
biholomorphic. Thus Theorem 1.C suggests a very interesting subject
to discover some interesting bounded weakly pseudoconvex domains
$D_1,\;D_2$ in $\mathbb{C}^n$ $(n\geq 2)$ such that any proper
holomorphic mapping from $D_1$ to $D_2$ is a biholomorphism. Even
though the bounded homogeneous domains in $\mathbb{C}^n$ are always
pseudoconvex, there are, of course, many such domains (e.g., all
bounded symmetric domains of rank $\geq 2$) without smooth boundary.
The lack of boundary regularity usually presents a serious
analytical difficulty. In 1984, by using the results of Bell
\cite{Bell} and Tumanov-Henkin \cite{Tuma}, Henkin and Novikov
\cite{Henkin} proved the following result (see Th.3.3 in
Forstneri\v{c} \cite{F} for references).

\vskip 6pt

\noindent {\bf Theorem 1.D} (Henkin and Novikov \cite{Henkin})  {\it
Any proper holomorphic self-mapping on an irreducible bounded
symmetric domain of rank $\geq 2$ is an analytic automorphism.}

\vskip 6pt

Further, using the idea in Mok-Tsai \cite{Mok-Tsai} and Tsai
\cite{Tsai}, Tu \cite{Tu1,Tu2} (one of the authors of the current
article) and Mok-Ng-Tu \cite{Mok-Ng} obtained some rigidity results
of proper holomorphic mappings between equidimensional bounded
symmetric domains (also called Cartan's domains). Recently,
Ahn-Byun-Park \cite{ABP} determined the automorphism group of the
Cartan-Hartogs domains (also called extended Cartan's domains) over
classical domains. In the past decade, Isaev \cite{Isa},
Isaev-Krantz \cite{IK} and Kim-Verdiani \cite{KV} also described the
automorphism groups of hyperbolic domains.

\vskip 6pt

The Fock-Bargmann-Hartogs domains $D_{n,m}(\mu)$ are defined by
$$D_{n,m}(\mu):=\{(z,w)\in\mathbf{C}^n\times\mathbf{C}^m:\|w\|^2<e^{-\mu\|z\|^2}\}, \;\;\; \mu>0.$$
The Fock-Bargmann-Hartogs domains $D_{n,m}(\mu)$ are unbounded
strongly pseudoconvex domains in $\mathbf{C}^{n+m}$. We note that
each $D_{n,m}(\mu)$ contains $\{(z, 0)\in
\mathbf{C}^n\times\mathbf{C}^m\} \cong \mathbf{C}^n$. Thus each
$D_{n,m}(\mu)$ is not hyperbolic in the sense of Kobayashi and
$D_{n,m}(\mu)$ can not be biholomorphic to any bounded domain in
$\mathbf{C}^{n+m}$. Therefore, each Fock-Bargmann-Hartogs domain
$D_{n,m}(\mu)$ is an unbounded non-hyperbolic domain in
$\mathbf{C}^{n+m}.$

\vskip 6pt

In 2013, Yamamori \cite{Yamamori} gave an explicit formula for the
Bergman kernel of the Fock-Bargmann-Hartogs domains in terms of the
polylogarithm functions. In 2014, by checking that the Bergman
kernel ensures revised Cartan's theorem, Kim-Ninh-Yamamori
\cite{Kim} determined the automorphism group of the
Fock-Bargmann-Hartogs domains as follows:

\vskip 6pt

\noindent {\bf Theorem 1.E} (Kim-Ninh-Yamamori \cite{Kim}) {\it The
automorphism group ${\rm Aut}(D_{n,m}(\mu))$ is exactly the group
generated by all automorphisms of  $D_{n,m}(\mu)$ as follows:
\begin{equation*}
\begin{array}{l}
\varphi_{U}:(z,w)\longmapsto(Uz,w),\quad U\in \mathcal{U}(n); \\\\
\varphi_{U^{'}}:(z,w)\longmapsto(z,U^{'}w),\quad U^{'}\in \mathcal{U}(m);\\\\
\varphi_{v}:(z,w)\longmapsto(z+v,e^{-\mu \langle z,v
\rangle-\frac{\mu}{2} {\left\lVert v\right\rVert}^{2}}w),\quad (v\in
\mathbb{C}^{n}),
\end{array}
\end{equation*}
where $\mathcal{U}(k)$ is the unitary group of degree $k,$ and
$\langle \cdot,\cdot \rangle$ is the standard Hermitian inner
product on $\mathbb{C}^{n}$.}

\vskip 6pt

The purpose of this article is to  prove the rigidity result on
proper holomorphic mappings between equidimensional
Fock-Bargmann-Hartogs domains as follows.

\begin{thm}
{If $D_{n,m}(\mu)$ and $D_{n',m'}(\mu')$ are two equidimensional
Fock-Bargmann-Hartogs domains with $m\geq 2$ and $f$ is a proper
holomorphic mapping of $D_{n,m}(\mu)$ into $D_{n',m'}(\mu')$, then
$f$ is a biholomorphism  between $D_{n,m}(\mu)$ and
$D_{n',m'}(\mu')$.}
\end{thm}

\noindent {\bf For example.}  Let
$\Phi(z_1,\cdots,z_n,w_1):=(\sqrt{2}z_1,\cdots,\sqrt{2}z_n,w_1^2)$,
$(z_1,\cdots,z_n,w_1)\in D_{n,1}(\mu)$. Then $\Phi$ is a proper
holomorphic self-mapping of $D_{n,1}(\mu)$, but it is  branched and
isn't an automorphism of $D_{n,1}(\mu)$. Then the assumption "$m\geq
2$" in Theorem 1.1 cannot be removed. Also, this example implies
that a proper holomorphic self-mapping of unbounded strongly
pseudoconvex domain in $\mathbf{C}^n(n\geq 2)$ is possibly not an
automorphism, and therefore, in general, Theorem 1.C does not hold
for unbounded strongly pseudoconvex domains in $\mathbf{C}^n(n\geq
2)$.

\vskip 6pt

Next we give a description of the biholomorphisms between two
Fock-Bargmann-domains as follows:

\begin{thm} {\it Let $D_{n,m}(\mu)$ and
$D_{n',m'}(\mu')$ be two equidimensional Fock-Bargmann-Hartogs
domains and let $f$ be a biholomorphism between $D_{n,m}(\mu)$ and
$D_{n',m'}(\mu')$. Then $n=n', \; m=m'$  and  therefore, there
exists $\varphi\in {\rm Aut}(D_{n',m'}(\mu'))$ such that
\begin{equation}\label{eq1}
f(z_1,\cdots,z_n,w_1,\cdots,w_m)=\varphi
(\sqrt{\mu/\mu'}z_1,\cdots,\sqrt{\mu/\mu'}z_n,w_1,\cdots,w_m).
\end{equation} }
\end{thm}

Now we shall present an outline of the argument in our proof of the
main results. Let  $f:\; D_{n,m}(\mu) \rightarrow D_{n',m'}(\mu')$
be  a proper holomorphic mapping between  two equidimensional
Fock-Bargmann-Hartogs domains. In order to prove that $f:\;
D_{n,m}(\mu) \rightarrow  D_{n',m'}(\mu')$ is a biholomorphism, it
suffices to show that $f$ is unbranched. Our proof consists of two
steps:

The first is to prove that $f$ extends holomorphically to their
closures. The transformation rule for Bergman kernels under proper
holomorphic mapping (e.g., Th. 1 in Bell \cite{Bell82}) is also
valid for unbounded domain (e.g., see Cor. 1 in Trybula
\cite{Trybula}). Note that the coordinate functions play a key role
in the approach of Bell \cite{Bell82} to extend proper holomorphic
mapping, but, in general, are no longer square integrable on
unbounded domains. In order to overcome the difficulty, by combining
the transformation rule for Bergman kernel under proper holomorphic
mapping in Bell \cite{Bell82} and an explicit form of the Bergman
kernel function for $D_{n,m}(\mu)$ in Yamamori \cite{Yamamori}, we
use a kind of semi-regularity at the boundary of the Bergman kernel
associated to $D_{n,m}(\mu)$ (see Th. 2.3 in this paper) to extend
the proper map holomorphically to a neighborhood of the closure
$\overline{D_{n,m}(\mu)}$ of $D_{n,m}(\mu)$, and then finish the
first step.

The second is to prove that $f:\; D_{n,m}(\mu) \rightarrow
D_{n',m'}(\mu')$ is unbranched assuming that the first step is
achieved. Assume that the zero locus  $S$ of the complex Jacobian of
the proper holomorphic mapping $f$ on $D_{n,m}(\mu)$ is not empty.
Then $S$ is of the codimension $1$. To finish the second step, by
using the strongly pseudoconvex boundary of $D_{n,m}(\mu)$ and the
local regularity for the mappings between strongly pseudoconvex
hypersurfaces (e.g., see Pin\v{c}uk \cite{Pin}), we get
$\overline{S}\cap b D_{n,m}(\mu)=\emptyset$ (note this will force
${S}$ to be compact if $D_{n,m}(\mu)$ is bounded) and then ${S}$ is
a complex analytic subset of $\mathbf{C}^{n+m}$. Further, we get
that the complex analytic subset $S$ of $\mathbf{C}^{n+m}$ must be
an algebraic set by its growth estimates. And, by considering the
dimension of the intersection of the projective closure
$\overline{S}$ of the affine algebraic set $S$ with the hyperplane
at infinity, we obtain that $S$ is of the codimension $\geq m$,
which forces $S$ to be $\emptyset$ by the assumption $m \geq 2$ in
Theorem 1.1. Therefore, $f$ is unbranched and is a biholomorphism.
This is the key ideas in proving the main results.

Our main work implies that any proper holomorphic self-mapping on
the Fock-Bargmann-Hartogs domain $D_{n,m}(\mu)$  with $m\geq 2$ must
be an automorphism.

\section{Preliminaries }
\subsection{  Bergman kernel associated to $D_{n,m}(\mu)$}

In this section we will make an investigation on a kind of
semi-regularity at the boundary of the Bergman kernel associated to
$D_{n,m}(\mu)$.

For a domain $\Omega$ in $\mathbf{C}^n$,  let $A^2(\Omega)$ be the
Hilbert space of square integrable holomorphic functions on $\Omega$
with the inner product:
$$\langle f,g\rangle=\int_{\Omega}f(z)\overline{g(z)} dV(z)\;\;(f,g\in \mathcal{O}(\Omega)),$$
where $dV$ is the Euclidean volume form. The Bergman kernel $K(z,w)$
of $A^2(\Omega)$ is defined as the reproducing kernel of the Hilbert
space $A^2(\Omega)$, that is, for all $f\in A^2(\Omega),$ we have
$$f(z)=\int_{\Omega}f(w)K(z,w)dV(w) \;\;(z\in\Omega).$$
For a positive continuous function $p$ on $\Omega$, let
$A^2(\Omega,p)$ be the weighted Hilbert space of square integrable
holomorphic functions with respect to the weight function $p$ with
the inner product:
$$\langle f,g\rangle=\int_{\Omega}f(z)\overline{g(z)}p(z)dV(z)\;\;\; (f,g\in \mathcal{O}(\Omega)).$$
Similarly, the weighted Bergman kernel $K_{A^2(\Omega,p)}$ of
$A^2(\Omega,p)$ is defined as the reproducing kernel of the Hilbert
space $A^2(\Omega,p)$. For a positive integer $m$, define the
Hartogs domain $\Omega_{m,p}$ over $\Omega$  by
$$\Omega_{m,p}=\{(z,w)\in\Omega\times\mathbf{C}^m: \|w\|^2<p(z)\}.$$

Ligocka \cite{Ligo} showed that the Bergman kernel of $\Omega_{m,p}$
can be expressed as infinite sum in terms of the weighted Bergman
kernel of $A^2(\Omega,p^k)\;(k=1,2,\cdots)$ as follows.

\vskip 6pt \noindent {\bf Theorem 2.1} (Ligocka \cite{Ligo}) {\it
Let $K_m$ be the Bergman kernel of $\Omega_{m,p}$ and let
$K_{A^2(\Omega,p^k)}$ be the weighted Bergman kernel of
$A^2(\Omega,p^k)$ $(k=1,2,\cdots)$. Then
$$K_m((z,w),(t,s))=\frac{m!}{\pi^m}\sum_{k=0}^{\infty}\frac{(m+1)_k}{k!}K_{A^2(\Omega,p^{k+m})}(z,t) \langle w,s\rangle^k,$$
where $(a)_k$ denotes the Pochhammer symbol
$(a)_k=a(a+1)\cdots(a+k-1)$. }

\vskip 6pt The Fock-Bargmann space is the weighted Hilbert space
$A^2(\mathbf{C}^n,e^{-\mu\|z\|^2})$ on $\mathbf{C}^n$ with the
Gaussian weight function $e^{-\mu\|z\|^2}\;(\mu>0)$. The reproducing
kernel of $A^2(\mathbf{C}^n,e^{-\mu\|z\|^2})$, called the
Fock-Bargmann kernel, is ${\mu^n e^{\mu \langle
z,t\rangle}}/{\pi^n}$ (see Bargmann \cite{Barg}). In 2013,  using
Th. 2.1 and the expression of the Fock-Bargmann kernel, Yamamori
\cite{Yamamori} give the Bergman kernel of the Fock-Bargmann-Hartogs
domain $D_{n,m}(\mu)$ as follows.

\vskip 6pt \noindent {\bf Theorem 2.2} (Yamamori \cite{Yamamori})
{\it The Bergman kernel of the Fock-Bargmann-Hartogs domain
$D_{n,m}(\mu)$ is given by
\begin{align*}
& \hskip 12pt K_{D_{n,m}(\mu)}((z,w),(t,s))\\
&=\frac{m!\mu^n}{\pi^{m+n}}\sum_{k=0}^{\infty}\frac{(m+1)_k(k+m)^n}{k!}
e^{\mu(k+m)\langle z,t\rangle}
\langle w,s\rangle^k\\
&=\frac{m!\mu^n}{\pi^{m+n}}\sum_{k=0}^{\infty}\frac{(m+1)_k(k+m)^n}{k!}
e^{\mu(k+m)\langle z,t\rangle} \sum_{\alpha_1+\cdots+\alpha_m=k}
\frac{k!}{\alpha_1!\cdots\alpha_m!}(w_1\overline{s_1})^{\alpha_1}\cdots (w_m\overline{s_m})^{\alpha_m}\\
&=\frac{m!\mu^n}{\pi^{m+n}}\sum_{\alpha\in\mathbf{N}^m}
\frac{(m+1)_{|\alpha|}(|\alpha|+m)^n}{\alpha!}
e^{\mu(|\alpha|+m)\langle z,t\rangle}
w^{\alpha}\overline{{s}^{\alpha}},
\end{align*}
where $(a)_k$ denotes the Pochhammer symbol
$(a)_k=a(a+1)\cdots(a+k-1)$ and
$\alpha=(\alpha_1,\cdots,\alpha_m)\in\mathbf{N}^m$ are multi-indices
of non-negative integer, $|\alpha|=\alpha_1+\cdots+\alpha_m$,
$\alpha!=\alpha_1!\cdots\alpha_m!$ and
$w^{\alpha}=w_1^{\alpha_1}\cdots w_m^{\alpha_m}$.}

\vskip 6pt Now we give a kind of semi-regularity at the boundary of
the Bergman kernel associated to $D_{n,m}(\mu)$ as follows:

\vskip 10pt \noindent {\bf Theorem 2.3} {\it Let $D_{n,m}(\mu)$ be a
Fock-Bargmann-Hartogs domain and let $K_{D_{n,m}(\mu)}((z,w),(t,s))$
be its Bergman kernel. If $E$ is a compact subset of $D_{n,m}(\mu)$,
then there is an open set $G$ containing $\overline{D_{n,m}(\mu)}$
such that for each $(t,s)\in E$, the function
$K_{D_{n,m}(\mu)}((z,w),(t,s))$ extends to be holomorphic on $G$ as
a function of $(z,w)$.}

\vskip 6pt \noindent {\bf Proof}. Since $E$ is a compact subset of
$D_{n,m}(\mu)$, there exists a real number $r$ with $0<r<1$ such
that $E\subset \{(z,w)\in D_{n,m}(\mu):\|w\|^2<r^2
e^{-\mu\|z\|^2}\}$. Let $G:=\{(z,w)\in
D_{n,m}(\mu):\|w\|^2<\frac{1}{r^2}e^{-\mu\|z\|^2}\}$. Then $G$ is an
open set containing $\overline{D_{n,m}(\mu)}$. By Theorem 2.2, we
have
$$K_{D_{n,m}(\mu)}((z,w),(t,s))=K_{D_{n,m}(\mu)}((z,rw),(t,\frac{1}{r}s)) $$
for all $(z,w)\in D_{n,m}(\mu)$, $(t,s)\in \{(z,w)\in
D_{n,m}(\mu):\|w\|^2<r^2 e^{-\mu \|z\|^2}\}$. Thus, for every fixed
$(t,s)\in E$, $K_{D_{n,m}(\mu)}((z,w),(t,s))$ extends
holomorphically to $G$ as a function of $(z,w)$. The proof of
Theorem 2.3 is finished.

\subsection{ Holomorphic extensions of proper holomorphic mappings }

In this section we will use Bell's transformation rule for Bergman
kernels under the proper holomorphic mappings and the
semi-regularity at the boundary of the Bergman kernel associated to
$D_{n,m}(\mu)$ to show that any proper holomorphic mapping $f$
between two equidimensional Fock-Bargmann-Hartogs domains
$D_{n,m}(\mu)$ and $D_{n',m'}(\mu')$ can be extended holomorphically
to the closure $\overline{D_{n,m}(\mu)}$ of $D_{n,m}(\mu)$.

The transformation rule for Bergman kernels under the proper
holomorphic mappings in Bell \cite{Bell82} plays an important role
in holomorphic extensions of proper holomorphic mappings. The
transformation rule (e.g., Th. 1 in Bell \cite{Bell82}) is also
valid for unbounded domain (e.g., see Cor. 1 in Trybula
\cite{Trybula}). Then we have the transformation rule for Bergman
kernels under the proper holomorphic mappings as follows.

\vskip 6pt \noindent {\bf Theorem 2.4} (Bell \cite{Bell82}, Theorem
1) {\it Suppose that $\Omega_1$ and $\Omega_2$ are two domains
$($not necessarily bounded$)$ in $\mathbf{C}^n$ and that $f$ is a
proper holomorphic mapping of $\Omega_1$ onto $\Omega_2$ of order
$r$. Let $u=\det [f']$ and let $F_1, F_2, \cdots, F_r$ denote the
$r$ local inverses to $f$ defined locally on $\Omega_2\setminus S$
where $S=\{f(z):u(z)=0\}$. Let $U_k=\det [F_k']$ and let $K_i(z,w)$
denote the Bergman kernel function associated to $\Omega_i$ for
$i=1,2$. The Bergman kernels transform according to
\begin{equation}\label{eq2}
\sum_{k=1}^rK_1 (z, F_k(w))\overline{U_k(w)}=u(z)K_2(f(z),w)
\end{equation}
for all $z\in \Omega_1$ and $w\in \Omega_2\setminus S$.}

\vskip 6pt

\noindent {\bf Remark on Theorem 2.4}. The removable singularity
theorem states that if $V$$(\varsubsetneqq D)$ is a complex variety
in a domain $D$ and $h\in L^2(D)$ (i.e., The Hilbert space of square
integrable functions on $D$) is holomorphic on $D\backslash V$, then
$h$ is holomorphic on $D$. Then the function on the left-hand side
of \eqref{eq2}  extends to be antiholomorphic in $w$ for all
$w\in\Omega_2$ by the removable singularity theorem (see Bell
\cite{Bell82} for references here).

\vskip 6pt

Now we will use Bell's transformation rule for Bergman kernels and
the semi-regularity at the boundary of the Bergman kernel associated
to $D_{n,m}(\mu)$ to show the holomorphic extension theorem as
follows.

\vskip 6pt \noindent {\bf Theorem 2.5} {\it If $D_{n,m}(\mu)$ and
$D_{n',m'}(\mu')$ are two equidimensional Fock-Bargmann-Hartogs
domains and $f$ is a proper holomorphic mapping of $D_{n,m}(\mu)$
into $D_{n',m'}(\mu')$, then $f$ extends to be holomorphic in a
neighborhood of $\overline{D_{n,m}(\mu)}$.}

\vskip 6pt \noindent {\bf Proof}. Let $f$ be a proper holomorphic
mapping of $D_{n,m}(\mu)$ onto $D_{n',m'}(\mu')$ and $u=\det[f']$. A
classical theorem due to R. Remmert (c.f. Rudin \cite{Rudin},
Theorem 15.1.9) states that $f$ is a branched covering of some
finite order $r$ and the set $S=\{f(z,w)\in D_{n',m'}(\mu'):
u(z,w)=0\}$ is a complex analytic variety in $D_{n',m'}(\mu')$.

Let $F_1,F_2,\cdots,F_r$ denote the $r$ local inverses to $f$
defined locally on $D_{n',m'}(\mu')\setminus S$. Let
$U_k=\det[F_k']$ and let $K_1((z,w),(t,s))$ and
$K_2((z',w'),(t',s'))$ denote the Bergman kernel function associated
to $D_{n,m}(\mu)$ and $D_{n',m'}(\mu')$ respectively. Write
$$H((z,w),(t',s'))=\sum_{k=1}^r
K_1((z,w),F_k(t',s'))\overline{U_k(t',s')}.$$ Then, from Remark on
Theorem 2.4, we have that $H((z,w),(t',s'))$ is holomorphic in
$(z,w)$  and is antiholomorphic in $(t',s')$ for all $((z,w),
(t',s'))\in  D_{n,m}(\mu)\times D_{n',m'}(\mu').$

With this notation, the transformation formula (\ref{eq2}) for
Bergman kernels becomes
\begin{equation}\label{eq5}
H((z,w),(t',s'))=u(z,w)K_2(f(z,w),(t',s')).
\end{equation}
Write
$f(z,w)=(f_1(z,w),f_2(z,w))\in\mathbf{C}^{n'}\times\mathbf{C}^{m'}.$
For $\alpha=(\alpha',\alpha'')\in\mathbf{N}^n\times \mathbf{N}^m$,
write
$$H^{(\alpha)}((z,w),(t',s')):=\frac{\partial^{\alpha}}{\partial
\overline{t'}^{\alpha'}\partial
\overline{s'}^{\alpha''}}H((z,w),(t',s')).$$  By differentiating the
equation (\ref{eq5}) with respect to
$(\overline{t'},\overline{s'})$, from Theorem 2.2, we have
\begin{align*}
&H^{(\alpha)}((z,w),(t',s'))=u(z,w) \frac{\partial^{\alpha}}{\partial \overline{t'}^{\alpha'}\partial \overline{s'}^{\alpha''}}K_2(f(z,w),(t',s'))\\
&=u(z,w) \frac{\partial^{\alpha}}{\partial \overline{t'}^{\alpha'}\partial \overline{s'}^{\alpha''}}
[\frac{m'!\mu'^{n'}}{\pi^{m'+n'}}\sum_{\beta\in \mathbf{N}^{m'}}\frac{(m'+1)_{|\beta|}(|\beta|+m')^{n'}}{\beta!}e^{\mu'(|\beta|+m')\langle f_1(z,w),t'\rangle}f_2(z,w)^{\beta}\overline{s'}^{\beta}]\\
&=u(z,w) \frac{m'!\mu'^{n'}}{\pi^{m'+n'}}\sum_{\beta:
\beta-\alpha''\in\mathbf{N}^{m'}}\frac{(m'+1)_{|\beta|}(|\beta|+m')^{n'}}{\beta!}\mu'^{|\alpha'|}(|\beta|+m')^{|\alpha'|}f_1(z,w)^{\alpha'}
e^{\mu'(|\beta|+m')\langle f_1(z,w),t'\rangle}\\
&\;\;\;\times \frac{\beta !}{(\beta-\alpha'')!} f_2(z,w)^{\beta}\overline{s'}^{\beta-\alpha''}.
\end{align*}
By putting $(t',s')=(0,0)$ in the above formula, we get
\begin{align*}
&H^{(\alpha)}((z,w),(0,0))\\
&=u(z,w)
\frac{m'!\mu'^{n'}}{\pi^{m'+n'}}\frac{(m'+1)_{|\alpha''|}(|\alpha''|+m')^{n'}}{\alpha''!}\mu'^{|\alpha'|}(|\alpha''|+m')^{|\alpha'|}f_1(z,w)^{\alpha'}
\alpha''! f_2(z,w)^{\alpha''}\\
&=u(z,w) \frac{m'!\mu'^{n'+|\alpha'|}}{\pi^{m'+n'}}(m'+1)_{|\alpha''|}(|\alpha''|+m')^{n'+|\alpha'|}f(z,w)^{\alpha}\\
&=C(\alpha)u(z,w)f(z,w)^{\alpha}
\end{align*}
for all
$\alpha=(\alpha',\alpha'')\in\mathbf{N}^{n'}\times\mathbf{N}^{m'}$.

Fix a neighborhood $V$ of $(0,0)$ with $V\subset\subset
D_{n',m'}(\mu')$. Then, for each $(t',s')\in\overline{V}\setminus
S$, we have $F_k(t',s')\in f^{-1}(\overline{V})\subset\subset
D_{n,m}(\mu)$ ($1\leq k\leq r$). Therefore, for all
$(t',s')\in\overline{V}\setminus S$, by Theorem 2.3, we have
$H((z,w),(t',s'))=\sum_{k=1}^r
K_1((z,w),F_k(t',s'))\overline{U_k(t',s')}$ can extends
holomorphically to a neighborhood $G$ of the closure
$\overline{D_{n,m}(\mu)}$ of $D_{n,m}(\mu)$ as a function of
$(z,w)$. Hence $H((z,w),(t',s'))$ is holomorphic in $(z,w)$ and
anti-holomorphic in $(t',s')$ for $((z,w),(t',s'))\in G\times
(V\setminus S)$.

Therefore, we have that $H((z,w),(t',s'))$ is holomorphic in $(z,w)$
and anti-holomorphic in $(t',s')$ for all $((z,w),(t',s'))\in
D_{n,m}(\mu)\times V$  (note $V \subset D_{n',m'}(\mu'))$  and for
all $((z,w),(t',s'))\in G\times (V\setminus S)$. So the Hartogs-type
extension theorem implies that $H((z,w),(t',s'))$ can be extended to
be a function on $G\times V$ which is holomorphic in $(z,w)$ and
anti-holomorphic in $(t',s')$ for all $((z,w),(t',s'))\in G\times
V$.

Hence $H^{(\alpha)}((z,w),(0,0))$ can extends holomorphically to the
neighborhood $G$ of $\overline{D_{n,m}(\mu)}$ of $D_{n,m}(\mu)$ as a
function of $(z,w)$ for all $\alpha\in\mathbf{N}^{n+m}$. Thus, the
function $u\cdot f^{\alpha}$ always extends holomorphically to the
neighborhood $G$ of $\overline{D_{n,m}(\mu)}$ for each
$\alpha\in\mathbf{N}^{n+m}$. This implies that $f$ extends to be
holomorphic in the neighborhood $G$ because the ring of germs of
holomorphic functions is a unique factorization domain. The proof of
Lemma 2.5 is finished.

\subsection{Cartan's theorem revisited}

Let $D\subset \mathbf{C}^N$ be a domain (not necessarily bounded)
with $0\in D$. Let $K_D(z,w)$ $(z,w\in D)$ be the Bergman kernel of
$D$. Let $T_D(z,w)$ be an $N\times N$ matrix defined by
\begin{eqnarray*}
\begin{aligned}
T_D(z,w):=
\begin{pmatrix}
\frac{\partial^2}{\partial\overline{w_1}\partial z_1}\log K_D(z,w) & \cdots & \frac{\partial^2}{\partial\overline{w_1}\partial z_N}\log K_D(z,w)  \\
\vdots & \ddots   & \vdots \\
\frac{\partial^2}{\partial\overline{w_N}\partial z_1}\log K_D(z,w) &
\cdots & \frac{\partial^2}{\partial\overline{w_N}\partial z_N}\log
K_D(z,w)
\end{pmatrix}.
\end{aligned}
\end{eqnarray*}
It is obviously that  $K_{D}(0,0)>0$ and $T_{D}(0,0)$ is positive
definite for any bounded domain $D\subset \mathbf{C}^N$.

 \vskip 6pt

Ishi-Kai \cite{Ishi} proved Cartan's theorem by using the notion of
the Bergman representative mapping for bounded circular domains.
However their proof is obviously applicable for an unbounded domain
whenever its Bergman kernel has some properties. Following the idea,
Kim-Ninh-Yamamori \cite{Kim} obtained a version of Cartan's theorem
for an unbounded circular domain $D\subset \mathbf{C}^N$ such that
$K_{D}(0,0)>0$ and $T_{D}(0,0)$ is positive definite, which assures
that any automorphism $f$ of  such an unbounded circular domain with
$f(0)=0$ must be linear. In this section we will get a slight
generalization of the result which states that any biholomorphism
$f$ between such two unbounded circular domains with $f(0)=0$ must
be linear.

 \vskip 6pt \noindent {\bf Lemma 2.6} (Ishi-Kai \cite{Ishi}, Prop.
2.1) {\it Let $D_{k}$ be a circular domain (not necessarily bounded)
in $\mathbf{C}^N$ with $0\in D_{k}\; (k=1,2).$ Let
$\varphi:D_1\rightarrow D_2 $ be a biholomorphism with
$\varphi(0)=0$. If $K_{D_k}(0,0)>0$ and $T_{D_k}(0,0)$ is positive
definite $(k=1,2)$, then $\varphi$ is linear.}

\vskip 6pt Remark. see Theorem 4 in Kim-Ninh-Yamamori \cite{Kim} for
references here.

\vskip 6pt

By the Lemmas 5 and 6 in Kim-Ninh-Yamamori \cite{Kim}, we have that
each Fock-Bargmann-Hartogs domain $D_{n,m}(\mu)$ satisfies the
conditions that $K_{D_{n,m}(\mu)}(0,0)>0$ and
$T_{D_{n,m}(\mu)}(0,0)$ is positive definite. Therefore, by Lemma
2.6,  we have a generalized Cartan's theorem for
Fock-Bargmann-Hartogs domains as follows:

\vskip 6pt

\noindent {\bf Theorem 2.7} {\it Let
$\varphi:D_{n,m}(\mu)\rightarrow D_{n',m'}(\mu')$ be a
biholomorphism between two equidimensional Fock-Bargmann-Hartogs
domains with $\varphi(0)=0$. Then $\varphi$ is linear. }

\subsection{Some lemmas about complex analytic sets }

In order to study the zero locus of the complex Jacobian of the
proper holomorphic mapping between two equidimensional
Fock-Bargmann-Hartogs domains, we need the following results.

\vskip 10pt \noindent {\bf Lemma 2.8} (Chirka \cite{Chirka}, \S 7.4
Theorem 3) {\it A pure $p$-dimensional analytic subset $A\subset
\mathbf{C}^n$ is algebraic if and only if it is contained, after
some unitary change of coordinates, in a domain
$D:\|z''\|<C(1+\|z'\|)^s$, where $z=(z',z''), z'=(z_1,\cdots,z_p)$,
and $C,s$ are certain constants.}

\vskip 10pt \noindent {\bf Lemma 2.9} (Chirka \cite{Chirka}, \S 7.2
Proposition 2) {\it The closure in $\mathbf{P}^n$ of an affine
algebraic set $A=\{\zeta\in\mathbf{C}^n:p(\zeta)=0\}$, where $p$ is
a polynomial of degree $s$, coincides with the projective algebraic
set $\{[z]\in\mathbf{P}^n:p^{*}(z)=0\}$, where $p^{*}$ is the
projectivization of $p$.}

\vskip 10pt In order to estimate the dimension of the zero locus of
the complex Jacobian of the proper holomorphic mapping between two
equidimensional Fock-Bargmann-Hartogs domains, we need the following
formula for the dimension of the intersection of two algebraic sets.

\vskip 10pt \noindent {\bf Lemma 2.10} (see Shafarevich
\cite{Shafa}) {\it Let $X,Y\subset \mathbf{P}^N$ be irreducible
quasiprojective varieties with $\dim X=n$ and $\dim Y=m$. Then any
(nonempty) component $Z$ of $X\cap Y$ has $\dim Z\geq n+m-N.$}

\vskip 10pt

In order to prove our main conclusion, we need the the preliminary lemma about regularity for
the mappings between strongly pseudoconvex hypersurfaces due to Pin\v{c}uk \cite{Pin} as follows.

\vskip 10pt \noindent {\bf Lemma 2.11} (Pin\v{c}uk \cite{Pin}, Lemma
1.3) {\it Let $D_1,D_2\subset \mathbf{C}^n$ be two domains, $p\in
bD_1$, and let $U$ be a neighborhood of $p$ in $\mathbf{C}^n$ such
that $U\cap \overline{D_1}$ is connected. Suppose that the mapping
$f=(f_1,\cdots,f_n):U\cap \overline{D_1}\rightarrow \mathbf{C}^n$ is
continuously differentiable in $U\cap \overline{D_1}$ and
holomorphic in $U\cap D_1$ with $f(U\cap bD_1)\subset bD_2$. Take a
domain $V\subset \mathbf{C}^n$ with $f(U\cap D_1)\subset V$. Suppose
that $U\cap bD_1$ and $U\cap bD_2$ are strongly pseudoconvex
hypersurfaces in $\mathbf{C}^n$. Then either $f$ is constant or the
Jacobian $J_f(z)=\det(\frac{\partial f_i}{\partial z_j})$ does not
vanish in $U\cap bD_1$.}

\section{Proof of main results }

\noindent
{\bf Proof of Theorem 1.1}.

Let $f:D_{n,m}(\mu)\rightarrow D_{n',m'}(\mu')$ be a proper
holomorphic mapping between two equidimensional
Fock-Bargmann-Hartogs domains $D_{n,m}(\mu)$ and $D_{n',m'}(\mu')$
with $m\geq 2$. By Theorem 2.5, $f$ extends holomorphically to a
neighborhood $V$ of $\overline{D_{n,m}(\mu)}$ with
\begin{equation}\label{eq6}
f(bD_{n,m}(\mu))\subset bD_{n',m'}(\mu').
\end{equation}
Define
$$A:=\{\zeta\in V: J_{f}(\zeta)=0\},$$
where $J_f(\zeta)=\det(\partial f_i/\partial \zeta_j)(\zeta)$ is the
complex Jacobian determinant of
$$f(\zeta):=(f_1(\zeta),\cdots,\\f_{n'+m'}(\zeta))\;\;(\zeta\in V).$$
Since $D_{n,m}(\mu)$ and $D_{n',m'}(\mu')$ are strongly pseudoconvex
domains, the Jacobian $J_f(\zeta)$ does not vanish on
$bD_{n,m}(\mu)$ by Lemma 2.11. Then we have $A\cap
bD_{n,m}(\mu)=\emptyset.$ Let $S:=A\cap {D_{n,m}(\mu)}.$ Therefore,
we have
\begin{equation}\label{eq7}
S\subset D_{n,m}(\mu),\;  \overline{S}\cap bD_{n,m}(\mu)=\emptyset.
\end{equation}

If $S\neq \emptyset$, by (\ref{eq7}), we can view $S$ as a complex
analytic set defined in $\mathbf{C}^{n+m}$. Moreover, for each $p\in
S (\subset D_{n,m}(\mu))$, we have
$$|w_m(p)|^2\leq\|w(p)\|^2<e^{-\mu\|z(p)\|^2}\leq 1\leq
1+\|(z,w')(p)\|,$$ where $w=(w',w_m).$ That is, we have
$$S\subset \{(z,w',w_m)\in \mathbf{C}^{n+m}:|w_m|<(1+\|(z,w')\|)^{1/2}\}.$$
Therefore,  by Lemma 2.8, we have that $S$ must be an algebraic set
of $\mathbf{C}^{n+m}$. Take an irreducible component $S'$ of $S$.
Now we consider the closure $\overline{S'}$ of $S'$ in
$\mathbf{P}^{n+m}$. By Lemma 2.9, $\overline{S'}$ is an projective
algebraic set and $\dim \overline{S'}=\dim S'=\dim S=n+m-1$.

Now we use Lemma 2.10 to give an upper bound $n$ for $\dim S'$ and
get a contradiction with  $\dim S'=n+m-1$. Let
$[\zeta,z_1,\cdots,z_n,w_1,\cdots,w_m]$ be the homogeneous
coordinate in $\mathbf{P}^{n+m}$ and embed $\mathbf{C}^{n+m}$ into
$\mathbf{P}^{n+m}$ as the affine piece
$U_0=\{[\zeta,z,w]\in\mathbf{P}^{n+m},\zeta\neq 0\}$ by
$$(z_1,\cdots,z_n,w_1,\cdots,w_m)\hookrightarrow
[1,z_1,\cdots,z_n,w_1,\cdots,w_m].$$ Then
\begin{equation}\label{eq8}
D_{n,m}(\mu)\cap U_0=\left\{[\zeta,z,w]\in\mathbf{P}^{n+m}:\;
\zeta\neq 0,
\frac{\|w\|^2}{|\zeta|^2}<e^{-\mu\frac{\|z\|^2}{|\zeta|^2}}\right\}.
\end{equation}
Let $H=\mathbf{P}^{n+m}\setminus \mathbf{C}^{n+m}$ be the hyperplane
at infinity, that is $H=\{\zeta=0\}\subset\mathbf{P}^{n+m}$.
Consider the affine piece
$U_1=\{[\zeta,z,w]\in\mathbf{P}^{n+m},z_1\neq 0\}$ of
$\mathbf{P}^{n+m}$ with affine coordinate
$(\xi,\lambda_2,\cdots,\lambda_n,\eta_1,\cdots,\eta_m)$. Then
$\xi=\frac{\zeta}{z_1},\lambda_2=\frac{z_2}{z_1},\cdots,\lambda_n=\frac{z_n}{z_1},\eta_1=\frac{w_1}{z_1},\cdots,\eta_m=\frac{w_m}{z_1}$.
Since
$\frac{\|w\|^2}{|\zeta|^2}=\frac{\|w\|^2}{|z_1|^2}\big|\frac{z_1}{\zeta}\big|^2=\frac{\|\eta\|^2}{|\xi|^2}$
and
$e^{-\mu\frac{\|z\|^2}{|\zeta|^2}}=e^{-\mu\frac{\|z\|^2}{|z_1|^2}\big|\frac{z_1}{\zeta}\big|^2}
=e^{-\mu\frac{(1+|\lambda_2|^2+\cdots+|\lambda_n|^2)}{|\xi|^2}}$, by
(\ref{eq8}), we have
\begin{equation}\label{eq9}
\begin{array}{l}
D_{n,m}(\mu)\cap U_1 \cap
U_0\\
=\left\{(\xi,\lambda_2,\cdots,\lambda_n,\eta_1,\cdots,\eta_m)\in\mathbf{C}^{n+m}:
|\eta_1|^2+\cdots+|\eta_m|^2<|\xi|^2e^{-\mu\frac{1+|\lambda_1|^2+\cdots+|\lambda_n|^2}{|\xi|^2}}\right\}.
\end{array}
\end{equation}
Let $S'_1=\overline{S'}\cap U_1$ be the affine piece of
$\overline{S'}$ in $U_1$ and let $H_1=H\cap U_1=\{\xi=0\}$ be the
affine piece of the projective hyperplane $H$ in $U_1$. For each
$p\in S'_1\cap H_1$, there exists a sequence of points
$\{p_k\}\subset \overline{S'}\cap((U_1\cap U_0)\setminus H_1)$ such
that $p_k\rightarrow p\; (k\rightarrow \infty)$. Since
$\{p_k\}\subset D_{n,m}(\mu)\cap U_1\cap U_0$, by (\ref{eq9}), we
have
\begin{equation}\label{eq10}
\|\eta(p_k)\|^2<|\xi(p_k)|^2e^{-\mu\frac{1+|\lambda_2(p_k)|^2+\cdots+|\lambda_n(p_k)|^2}{|\xi(p_k)|^2}}.
\end{equation}
Since $p\in H$, we have $\xi(p)=0$ and $\xi(p_k)\rightarrow 0\;
(k\rightarrow \infty)$. Let $k\rightarrow \infty$ in (\ref{eq10}),
we get $\|\eta(p)\|^2=0$. Therefore, $S'_1\cap H_1\subset
\{\xi=0,\eta_1=\cdots=\eta_m=0\}$. Hence, $\dim(S'_1\cap H_1)\leq
n-1$.

Further, by Lemma 2.10, we have $$n-1 \geq \dim(S'_1\cap H_1)\geq
\dim S'_1+\dim H_1-(n+m)=\dim S'_1-1.$$ Thus, $\dim S'_1\leq n$.

Therefore, $n+m-1=\dim S=\dim S'=\dim\overline{S'}=\dim S'_1\leq n$.
Hence, $m\leq 1$, this is a contradiction with the assumption $m\geq
2$ of Theorem 1.1. This means $S=\emptyset$.

Thus $f:D_{n,m}(\mu)\rightarrow D_{n',m'}(\mu')$ is unbranched. Since each Fock-Bargmann-Hartogs domain is simply connected,
we get that $f:D_{n,m}(\mu)\rightarrow D_{n',m'}(\mu')$ is a biholomorphism. The proof of Theorem 1.1 is finished.

\vskip 6pt
\noindent
{\bf Proof of Theorem 1.2}.

For the completeness, here we will not assume Theorem 1.E to prove
Theorem 1.2. Our proof of Theorem 1.2 is divided as two steps:

{\bf Step 1}. Let $f:D_{n,m}(\mu)\rightarrow D_{n',m'}(\mu')$ be a
biholomorphical mapping between two equidimensional
Fock-Bargmann-Hartogs domains $D_{n,m}(\mu)$ and $D_{n',m'}(\mu').$
We will show that $n=n', m=m'$ and there exists a $\varphi'\in {\rm
Aut}(D_{n',m'}(\mu'))$ such that $\varphi'\circ f$ is a linear
isomorphism.

Let $\mathcal{V}=\{(z,0)\in\mathbf{C}^n\times\mathbf{C}^m\}\subset
D_{n,m}(\mu)$ and
$\mathcal{V'}=\{(z',0)\in\mathbf{C}^{n'}\times\mathbf{C}^{m'}\}\subset
D_{n',m'}(\mu')$. Put $f(z,0)=(g(z),h(z))$ and $h(z)=(h_1(z),\cdots,
h_{m'}(z))$. Then we have
\begin{equation}\label{eq20}
\sum_{i=1}^{m'}|h_i(z)|^2=\|h(z)\|^2<e^{-\mu'\|g(z)\|^2}\leq 1.
\end{equation}
It follows that $h_i$ is a bounded holomorphic function on
$\mathbf{C}^n$ for all $1\leq i\leq m'$. Then Liouville's theorem
implies that $h_i$ is constant. Since $g$ is a non-constant entire
function, $g$ is unbounded. Therefore, by \eqref{eq20},  $h$ must be
identically equal to zero. This means
$f(\mathcal{V})\subset\mathcal{V'}$. In a similar way we have
$f^{-1}(\mathcal{V'})\subset\mathcal{V}$. Thus, $f|_{\mathcal{V}}$
is a biholomorphism between $\mathcal{V}$ and $\mathcal{V'}$.
Therefore, we have $n=n',\; m=m'$.

Let $f(0,0)=(v,\tilde{v})$. Then $\tilde{v}=0$. Let $\varphi_{-v}\in
{{\rm Aut}}(D_{n',m'}(\mu'))$ be defined by
$$\varphi_{-v}(z',w')=(z'-v,e^{-\mu'\langle
z',-v\rangle-\frac{\mu}{2}\|-v\|^2}w').$$ Then $\varphi_{-v}\circ
f(0,0)=(0,0).$ By Theorem 2.7, $\varphi_{-v}\circ f$ is a linear
mapping.

{\bf Step 2}. We now prove that there exist $\varphi^*\in {\rm
Aut}(D_{n',m'}(\mu'))$ such that $\varphi^*\circ\varphi_{-v}\circ f$
has the desired form (\ref{eq1}).

Since $\varphi_{-v}\circ f$ is linear, the map $\varphi_{-v}\circ f$
can be written as a matrix form, namely,
\begin{eqnarray*}
\begin{aligned}
\varphi_{-v}\circ f(z,w)=
\begin{pmatrix}
A & C\\
D & B\\
\end{pmatrix}
\begin{pmatrix}
z\\w
\end{pmatrix},
\end{aligned}
\end{eqnarray*}
where $A\in M_{n\times n}(\mathbf{C}), B\in M_{m\times
m}(\mathbf{C}), C\in M_{n\times m}(\mathbf{C})$ and $D\in M_{m\times
n}(\mathbf{C})$. Since $\varphi_{-v}\circ f(\mathcal{V})\subset
\mathcal{V'}$, we have
$$D=\mathbf{0}\in M_{n\times m}(\mathbf{C});\;\;\; \det A\neq 0;\;\;\; \det B\neq 0.$$
Then it follows that $\varphi_{-v}\circ f(z,w)=(Az+Cw,Bw)$ for all $(z,w)\in D_{n,m}(\mu)$, and
$(\varphi_{-v}\circ f)^{-1}(z',w')=(A^{-1}z'-A^{-1}CB^{-1}w',B^{-1}w')$ for all $(z',w')\in D_{n',m'}(\mu')$.

Now we prove that $B\in \mathcal{U}(m)$. For $\|w\|<1$, we have
$(0,w)\in D_{n,m}(\mu)$ and $\varphi_{-v}\circ f(0,w)=(Cw,Bw)\in
D_{n',m'}(\mu')$. Thus $\|Bw\|<1$. On the other hand, for
$\|w'\|<1$, we have $(0,w')\in D_{n',m'}(\mu')$ and
$(\varphi_{-v}\circ f)^{-1}(0,w')=(-A^{-1}CB^{-1}w',B^{-1}w')\in
D_{n,m}(\mu)$. Thus $\|B^{-1}w'\|<1$. Therefore,
$$B:\mathbf{B}^m\rightarrow \mathbf{B}^m,\; w\rightarrow Bw $$
is a linear automorphism of the unit ball $\mathbf{B}^m$ in
$\mathbf{C}^m$, and thus $B\in \mathcal{U}(m)$.

Next we prove that $C=\mathbf{0}\in M_{n\times m}(\mathbf{C})$. For
$\|w\|=1$, we have $(0,w)\in bD_{n,m}(\mu)$. Since
$\varphi_{-v}\circ f$ is a linear biholomorphism, we have
$\varphi_{-v}\circ f(b D_{n,m}(\mu))=b D_{n',m'}(\mu')$. Hence,
$\varphi_{-v}\circ f(0,w)=(Cw,Bw)\in bD_{n',m'}(\mu')$. Therefore,
by $B\in \mathcal{U}(m)$, we have
$1=\|w\|^2=\|Bw\|^2=e^{-\mu'\|Cw\|^2}$. Thus, $Cw=0$ for all
$\|w\|=1$. Hence, $C=\mathbf{0}\in M_{n\times m}(\mathbf{C})$.

To complete our proof, it suffices to show that $A=\sqrt{\mu/\mu'}U$ for some $U\in \mathcal{U}(n)$.
For any $z\in\mathbf{C}^n$, we can take $w\in\mathbf{C}^m$ such that $(z,w)\in bD_{n,m}(\mu)$.
Since $\varphi_{-v}\circ f(bD_{n,m}(\mu))= bD_{n',m'}(\mu')$ and $C=\mathbf{0}$,
we have $\varphi_{-v}\circ f(z,w)=(Az,Bw)\in bD_{n',m'}(\mu')$.
Thus, by $B\in \mathcal{U}(m)$, we get $e^{-\mu'\|Az\|^2}=\|Bw\|^2=\|w\|^2=e^{-\mu\|z\|^2}$. Therefore,
$\|\sqrt{\mu'/\mu}Az\|=\|z\|$ for all $z\in\mathbf{C}^n$. Hence, $U=\sqrt{\mu'/\mu}A\in \mathcal{U}(n)$
and $A=\sqrt{\mu/\mu'}U$.

Let $\varphi_{U^{-1}},\; \varphi_{B^{-1}}\in {\rm
Aut}(D_{n',m'}(\mu'))$ be defined by
$$\varphi_{U^{-1}}(z',w')=(U^{-1}z',w'),\; \varphi_{B^{-1}}(z',w')=(z',B^{-1}w').$$
Then $\varphi:=\varphi_{U^{-1}}\circ \varphi_{B^{-1}}\circ
\varphi_{-v}\in {\rm Aut}(D_{n',m'}(\mu'))$, and we have
$$\varphi\circ f(z_1,\cdots,z_n,w_1,\cdots,w_m)=(\sqrt{\mu/\mu'}z_1,\cdots,\sqrt{\mu/\mu'}z_n,w_1,\cdots,w_m).$$
The proof of Theorem 1.2 is finished.

\vskip 10pt

\noindent\textbf{Acknowledgments}\quad The authors would like to
thank Professors Xiaojun Huang and Ngaiming Mok for helpful
suggestions. In addition, the authors are grateful to the referees
for many helpful comments. The project was supported by the National
Natural Science Foundation of China (No.11271291).

\end{document}